\documentclass{article}

\usepackage{graphicx} 
\usepackage{amssymb,amsmath,amsfonts,amsthm,lineno}
\usepackage[mathscr]{euscript}
\usepackage{latexsym}
\usepackage{graphics}
\usepackage{indentfirst}
\usepackage[colorlinks=true,allcolors=blue]{hyperref}

\usepackage{comment}
\usepackage[shortlabels]{enumitem}
\usepackage[english]{babel}
\usepackage{etoolbox}

\usepackage{esdiff}
\usepackage{thmtools}
\usepackage{thm-restate}

\usepackage{tikz}
\usetikzlibrary{fit,arrows.meta,backgrounds,shapes,decorations, graphs, graphs.standard,quotes, decorations.markings, arrows.meta, backgrounds}
\usepackage{placeins}

\newif\iflatexml
\newif\iffirstzref 
\latexmlfalse
\ifdefined\LaTeXML\latexmltrue\fi
\ifdefined\LaTeXMLversion\latexmltrue\fi

\iflatexml
  \makeatletter
  \newcommand\cref[1]{%
    \firstzreftrue
    \@for\@tmpc:=#1\do{%
      \iffirstzref\firstzreffalse\else, \fi
      \autoref{\@tmpc}}}
  \makeatother
  \let\Cref\cref
  \let\zcref\cref
\else
  \usepackage{zref-clever}
  \zcsetup{cap}
  \zcsetup{nameinlink}
  \zcRefTypeSetup{theorem}{Name-sg=Theorem, name-sg=Theorem, Name-pl=Theorems, name-pl=Theorems}
  \zcRefTypeSetup{lemma}{Name-sg=Lemma, name-sg=Lemma, Name-pl=Lemmas, name-pl=Lemmas}
  \zcRefTypeSetup{obs}{Name-sg=Observation, name-sg=Observation, Name-pl=Observations, name-pl=Observations}
  \zcRefTypeSetup{pro}{Name-sg=Proposition, name-sg=Proposition, Name-pl=Propositions, name-pl=Propositions}
  \zcRefTypeSetup{defn}{Name-sg=Definition, name-sg=Definition, Name-pl=Definitions, name-pl=Definitions}
  \zcRefTypeSetup{cor}{Name-sg=Corollary, name-sg=Corollary, Name-pl=Corollaries, name-pl=Corollaries}
  \zcRefTypeSetup{conj}{Name-sg=Conjecture, name-sg=Conjecture, Name-pl=Conjectures, name-pl=Conjectures}
  \zcRefTypeSetup{cl}{Name-sg=Claim, name-sg=Claim, Name-pl=Claims, name-pl=Claims}
  \zcRefTypeSetup{example}{Name-sg=Example, name-sg=Example, Name-pl=Examples, name-pl=Examples}
  \zcRefTypeSetup{quest}{Name-sg=Question, name-sg=Question, Name-pl=Questions, name-pl=Questions}
  \let\cref\zcref
\fi

\newtheorem*{theorem*}{Theorem}
\newtheorem{theorem}{Theorem}
\newtheorem{lemma}[theorem]{Lemma}
\newtheorem{pro}[theorem]{Proposition}

\newtheorem{cor}[theorem]{Corollary}
\newtheorem{conj}[theorem]{Conjecture}

\theoremstyle{definition}
\newtheorem{example}[theorem]{Example}

\newtheorem{innercustomlem}{Lemma}
\newenvironment{customlem}[1]
  {\renewcommand\theinnercustomlem{#1}\innercustomlem}
  {\endinnercustomlem}

\makeatletter
\newcommand{\forceitemlabel}[1]{%
  \refstepcounter{enumi}%
  \protected@write\@auxout{}{\string\newlabel{#1}{{\theenumi}{\thepage}}}%
}
\makeatother

\newif\ifshowfull

\makeatletter
\iflatexml
  \newlist{statement}{enumerate}{1}
  \setlist[statement]{
    label=\textup{(\alph*)},
    ref=\textup{(\alph*)},
    nosep
  }

  \newcommand\zrefparent[2]{%
    \@for\@tmpb:=#2\do{\expandafter\gdef\csname zrefparent@\@tmpb\endcsname{#1}}}

  \newcommand\zrefWithTheorem[1]{%
    \firstzreftrue
    \@for\@tmpa:=#1\do{%
      \iffirstzref\firstzreffalse\else, \fi
      \ifcsname zrefparent@\@tmpa\endcsname
        \expandafter\autoref\expandafter{\csname zrefparent@\@tmpa\endcsname}\ref{\@tmpa}%
      \else
        \ref{\@tmpa}%
      \fi}}

  \zrefparent{thm:examples}{examplesa,examplesb,examplesc,examplesf,examplesd,examplese}
  \zrefparent{lem:lboundlistcolorfunction}{bound1intro,bound2intro}
  \zrefparent{lem:mult2m-t}{statement1,statement2,statement3}
\else
  \newcommand*{\currenttheoremenv}{}

  \AtBeginEnvironment{theorem}{\xdef\currenttheoremenv{theorem}}
  \AtBeginEnvironment{lemma}{\xdef\currenttheoremenv{lemma}}
  \AtBeginEnvironment{proposition}{\xdef\currenttheoremenv{proposition}}
  \AtBeginEnvironment{corollary}{\xdef\currenttheoremenv{corollary}}

  \newcommand*\doublelabel[1]{\protect\@twolabels{#1}{\@currentlabel#1}}
  \let\@twolabels\@firstoftwo

  \newcommand*\versionWithTheorem[1]{%
    \@ifstar{\versionWithTheorem@aux{#1* me}}{\versionWithTheorem@aux{#1}}%
  }
  \newcommand*\versionWithTheorem@aux[2]{%
    \begingroup\let\@twolabels\@secondoftwo#1{#2}\endgroup%
  }

  \DeclareRobustCommand*\zrefWithTheorem{\versionWithTheorem\zcref}

  \newlist{statement}{enumerate}{1}
  \setlist[statement]{
    label=\textup{(\alph*)},
    ref=\doublelabel{\textup{(\alph*)}},
    before={\zcsetup{countertype={statementi=\currenttheoremenv}}},
    nosep
  }
\fi
\makeatother

\newcommand{\M}{\mathcal{M}}
\newcommand{\N}{\mathbb{N}}

\newcommand\C[1]{\mathscr{C}_{(#1)}}

\newcommand\K[1]{K_{2*#1}}

\newcommand\abs[1]{\left\lvert #1 \right\rvert}

\newcommand\ceil[1]{\left\lceil #1 \right\rceil}

\allowdisplaybreaks

\title{Enumerative Chromatic Choosability of Complete Multipartite Graphs with Small Partite Sets}

\author{Sarah Allred \thanks{Department of Mathematics and Statistics, University of South Alabama, Mobile, AL, USA (sarahallred@southalabama.edu)} 
\and Joseph T. Brown \thanks{School of Computing, University of South Alabama, Mobile, AL, USA (jtb2421@jagmail.southalabama.edu)}
\and Devin Dickinson \thanks{Department of Mathematics and Statistics, University of South Alabama, Mobile, AL, USA (dld1721@jagmail.southalabama.edu)} \and H. C. Fritz
\thanks{Department of Mathematics and Statistics, University of South Alabama, Mobile, AL, USA (hcf2021@jagmail.southalabama.edu)} \and Aiden McCain \thanks{Department of Mathematics and Statistics, University of South Alabama, Mobile, AL, USA (cam2326@jagmail.southalabama.edu)}\and Jeffrey A. Mudrock \thanks{Department of Mathematics and Statistics, University of South Alabama, Mobile, AL, USA (mudrock@southalabama.edu)}\and 
Daniel Smith \thanks{Department of Mathematics and Statistics, University of South Alabama, Mobile, AL, USA (dps2321@jagmail.southalabama.edu)}
\and Morgan Taylor \thanks{Department of Mathematics and Statistics, University of South Alabama, Mobile, AL, USA (mot2221@jagmail.southalabama.edu)} \and Benjamin Whatley
\thanks{Mitchell College of Business, University of South Alabama, Mobile, AL, USA (bvw2422@jagmail.southalabama.edu)}}
\date{}

\begin{document}

\maketitle

\begin{abstract}
Chromatic choosability is a notion of fundamental importance in list coloring.  A graph is \emph{chromatic-choosable} when its chromatic number is equal to its list chromatic number.  In their 1979 paper on list coloring, Erd\H{o}s, Rubin, and Taylor showed that complete $m$-partite graphs with every partite set of size 2, denoted by $K_{2*m}$, are chromatic-choosable.  In 1990, Kostochka and Sidorenko introduced the list color function of a graph $G$, denoted by $P_{\ell}(G,m)$, which is the list analogue of the chromatic polynomial of $G$, $P(G,m)$.  A graph $G$ is \emph{enumeratively chromatic-choosable} when $P_{\ell}(G,m) = P(G,m)$ whenever $m \in \N$.  In this paper, we explore the question as to whether $K_{2*m}$ satisfies this stronger property. We establish a recursive lower bound on $P_{\ell}(K_{2*m},m)$ and prove that $K_{2*3}$ and $K_{2*4}$ are \emph{weakly enumeratively chromatic-choosable}, meaning that equality holds at the chromatic number. We also develop a general lower bound for the list color function of certain joins with complete graphs and use it to prove that $K_p\vee K_{2*5}$ is weakly enumeratively chromatic-choosable whenever $p\geq 15$, and that $K_p\vee K_{2,2,3}$ is weakly enumeratively chromatic-choosable whenever $p\geq 3$. These results provide progress toward a weaker form of a conjecture of Allred and Mudrock concerning the enumerative chromatic choosability of complete graph joins.

\end{abstract}

\section{Introduction}\label{sec:Introduction}
In this paper, all graphs are finite, nonempty, and simple.  The set of natural numbers is $\N = \{1,2,3, \ldots \}$.  For $m \in \N$, we write $[m]$ for the set $\{1, \ldots, m \}$ and take $[0]$ to be the empty set.   For $m, n \in \N$ and $m \leq n$, we write $[m : n]$ for the set $\{m, \ldots, n\}$. If $m > n$, then $[m : n] = \emptyset$.  We use $K_{n*m}$ for the complete $m$-partite graphs where each partite set has size $n$. If $G$ and $H$ are vertex disjoint graphs, we write $G \vee H$ for the join of $G$ and $H$. Terms and notation that are not defined here follow West~\cite{W01}.

\subsection{Classical Coloring and List Coloring}\label{subsec:listcol}
In classical vertex coloring, a \emph{proper m-coloring} of a graph $G$ assigns a color to each $v \in V(G)$ and uses up to $m$ colors from $[m]$ such that adjacent vertices receive different colors.  The \emph{chromatic number} of $G$, denoted by $\chi(G)$, is the smallest $m$ such that $G$ has a proper $m$-coloring.
List coloring is a well-known variation on classical vertex coloring which was introduced independently by Vizing~\cite{V76} and Erd\H{o}s, Rubin, and Taylor~\cite{ET79} in the 1970s.  A \emph{list assignment} $L$ for a graph $G$ assigns to each vertex $v \in V(G)$ a list of available colors $L(v)$.  We say $G$ is \emph{$L$-colorable} if there is a proper coloring $f$ of $G$ such that $f(v) \in L(v)$ for each $v \in V(G)$ (we refer to $f$ as a \emph{proper $L$-coloring} of $G$).  A list assignment $L$ for $G$ is called an \emph{$m$-assignment} if $|L(v)|=m$ for each $v \in V(G)$.  A graph $G$ is said to be $m$-choosable if it is $L$-colorable whenever $L$ is an $m$-assignment for $G$.  The \emph{list chromatic number} of a graph $G$, denoted by $\chi_\ell(G)$, is the smallest $k$ for which $G$ is $k$-choosable.  It is easy to show that for any graph $G$, $\chi(G) \leq \chi_\ell(G)$.  Moreover, the gap between the chromatic number and list chromatic number of a graph can be arbitrarily large since $\chi_{\ell}(K_{n,t}) = n+1$ whenever $t \geq n^n$ (see e.g.,~\cite{M18} for further details). 

A graph $G$ is called \emph{chromatic-choosable} if $\chi_{\ell}(G) = \chi(G)$~\cite{O02}.   Determining whether a graph is chromatic-choosable is, in general, a challenging problem.  In their 1979 paper introducing list coloring, Erd\H{o}s, Rubin, and Taylor~\cite{ET79} showed that $K_{2*m}$ is chromatic-choosable.  In addition, perhaps the most well-known conjecture involving list coloring concerns chromatic choosability.  Indeed, the famous Edge List Coloring Conjecture states that every line graph of a loopless multigraph is chromatic-choosable (see~\cite{HC92}).

While the chromatic and list chromatic numbers concern the existence of proper colorings, we may also ask how many proper $m$-colorings or proper $L$-colorings a graph admits. The chromatic polynomial addresses the former question. For each $m \in \N$, the \emph{chromatic polynomial} of a graph $G$, denoted by $P(G,m)$, is the number of proper $m$-colorings of $G$. The function $P(G,m)$ is a polynomial in $m$ of degree $|V(G)|$ (e.g., see~\cite{DKT05}). For example, whenever $n \in \N$ it is well-known that $P(K_n,m) = \prod_{i=0}^{n-1} (m-i)$, $P(T,m) = m(m-1)^{n-1}$ whenever $T$ is a tree on $n$ vertices, and $P(K_1 \vee G, m) = m P(G,m-1)$ for every graph $G$ and each $m \geq 2$ (see~\cite{B94} and~\cite{W01}). Recall that a graph $G$ with $\chi(G)=m$ is \emph{uniquely $m$-colorable} if there is a unique partition of $V(G)$ into $m$ independent sets. An important fact for this paper is that if $G$ is uniquely $m$-colorable, then $P(G,m)=m!$.

In the early 1990s, Kostochka and Sidorenko~\cite{AS90} introduced an analogue of the chromatic polynomial for list coloring. If $L$ is a list assignment for $G$, let $P(G,L)$ denote the number of proper $L$-colorings of $G$. The \emph{list color function} $P_\ell(G,m)$ is the minimum value of $P(G,L)$ where the minimum is taken over all possible $m$-assignments $L$ for $G$.  Since an $m$-assignment could assign the same $m$ colors to every vertex in a graph, it follows that $P_\ell(G,m) \leq P(G,m)$ for each $m \in \N$.  In general, the list color function can differ significantly from the chromatic polynomial for small values of $m$.  For example, for any $n \geq 2$, $P_{\ell}(K_{n,n^n},m) = 0$ and $P(K_{n,n^n},m) > 1$ whenever $m \in [2:n]$.  On the other hand, Donner~\cite{D92} showed in 1992 that for every graph $G$ there is a $k \in \N$ such that $P_\ell(G,m) = P(G,m)$ whenever $m \geq k$.

\subsection{Enumerative Chromatic Choosability}\label{subsec:ecc}
We now consider an enumerative analogue of chromatic choosability. A graph $G$ is said to be \emph{enumeratively chromatic-choosable} if $P_\ell(G,m)=P(G,m)$ for each $m\geq\chi(G)$. Since both functions vanish whenever $m<\chi(G)$, this is equivalent to requiring equality for every $m\in\mathbb{N}$.

It was recently proven by Zhang and Dong~\cite{DZ26} that for each $k \geq 3$ there is a graph $G$ with $\chi(G) = k$ such that $P_{\ell}(G,k) = P(G,k)$ and $P_{\ell}(G,k+1) < P(G, k+1)$.  So, if one proves $P_{\ell}(G, \chi(G)) = P(G, \chi(G))$ for some $G$, it does not imply that $G$ is enumeratively chromatic-choosable. 
Prior to this result, Allred and Mudrock \cite{AM25} introduced the notion of \emph{weak enumerative chromatic choosability}.  A graph $G$ is weakly enumeratively chromatic-choosable if $P_{\ell}(G, \chi(G)) = P(G, \chi(G))$. By Zhang and Dong's result, there are graphs that are weakly enumeratively chromatic-choosable but not enumeratively chromatic-choosable.

Relatively few classes of graphs are known to be enumeratively chromatic-choosable. In order to present some known examples, we need a definition. A \emph{theta graph}, denoted by $\Theta(l_1,l_2,l_3)$, consists of a pair of end vertices joined by three internally disjoint paths of lengths $l_1,l_2,l_3 \in \N$, where we allow at most one path to have length one.

The following theorem gives some known results on enumerative chromatic choosability.
\begin{theorem} \label{thm:examples} \cite{AM25, BK23, CL26, KK23, KM18, KN16, AS90}
The following statements hold.
\begin{statement}
    \item\label{examplesa} Chordal graphs are enumeratively chromatic-choosable.
    \item\label{examplesb} Cycles are enumeratively chromatic-choosable.
    \item\label{examplesc} $K_{2,3}$ is enumeratively chromatic-choosable.
     \item\label{examplesf} For every $k \in \N$, $K_1\vee \Theta(2,2,2k)$ is enumeratively chromatic-choosable.
    \item\label{examplesd} Suppose $G=\Theta(l_1,l_2,l_3)$ where $\min\{l_1,l_2,l_3\} = l_1$, $l_2 \geq 2$, and $l_3 \geq 2$.  Then $G$ is not enumeratively chromatic-choosable if and only if $l_1,l_2,$ and $l_3$ all have the same parity and $(l_1,l_2,l_3) \neq (2,2,2)$.
    \item\label{examplese} If $G$ is enumeratively chromatic-choosable (resp. weakly enumeratively chromatic-choosable), then $G \vee K_n$ is enumeratively chromatic-choosable (resp. weakly enumeratively chromatic-choosable) for any $n \in \N$. 
\end{statement}
\end{theorem}

The motivation for this paper is twofold. First, inspired by the result of Erd\H{o}s, Rubin, and Taylor~\cite{ET79}, we wished to determine whether $K_{2*m}$ is enumeratively chromatic-choosable or, at least, weakly enumeratively chromatic-choosable.  Notice that \zrefWithTheorem{examplesa,examplesb} imply $K_{2*m}$ is enumeratively chromatic-choosable when $m \in [2]$.  Second, Ohba~\cite{O02} showed that for any graph $G$, there exists an $N \in \N$ such that $K_p \vee G$ is chromatic-choosable whenever $p \geq N$.  Motivated by this, Allred and Mudrock made the following conjecture.

\begin{conj} \label{conj:join} \cite{AM25}
For any graph $G$, there is a $p \in \N$ such that $K_p \vee G$ is enumeratively chromatic-choosable.    
\end{conj}

With \cref{conj:join} in mind, we investigate, for a chromatic-choosable graph $G$ that is not known to be weakly enumeratively chromatic-choosable, whether there exists a $p \in \N$ such that $K_p \vee G$ is weakly enumeratively chromatic-choosable.

\subsection{Outline of the Paper}\label{subsec:outline}
In \cref{sec:prelimresults}, we begin by introducing tools we will use for the rest of the paper.  Then, we pursue a direction inspired by \cref{conj:join}.  Specifically, we show that when $G$ is a uniquely $m$-colorable graph satisfying certain conditions and not known to be weakly enumeratively chromatic-choosable, a lower bound on $P_{\ell}(K_p \vee G,m+p)$ must hold.  Indeed, we prove the following.
\begin{restatable}{theorem}{thmjoincomplete}\label{thm:joincomplete}
Suppose $G$ is a uniquely $m$-colorable graph where $m \geq 2$ and the unique partition of $V(G)$ into $m$ independent sets is $\{A_i : i \in [m]\}$.  
Suppose $G$ is such that $G - \bigcup_{i=1}^{m-1} A_{x_i}$ is weakly enumeratively chromatic-choosable whenever $(x_1, \ldots, x_{m-1}) \in [m]^{m-1}$ and $P_{\ell}(G,m) \geq a_0$.  For each $p \in \N$ let
\[a_p = \left( m \left( \frac{m+p}{m+p-1}\right) +p \right)a_{p-1}.\]
Then $P_{\ell}(K_p \vee G, m+p) \geq \min\{a_p, (m+p)!\}$.
\end{restatable}

Notice that the entries of $(x_1,\ldots,x_{m-1})$ in \cref{thm:joincomplete} need not be distinct. Consequently, the hypothesis requires that $G-\bigcup_{i\in I}A_i$ is weakly enumeratively chromatic-choosable for every nonempty proper subset $I$ of $[m]$.

Notice that $K_p \vee G$ is uniquely $(m+p)$-colorable which means that if $a_p \geq (m+p)!$, then $K_p \vee G$ is weakly enumeratively chromatic-choosable. \cref{cor:joincomplete} gives a condition on $a_0$ and $m$ that guarantees $a_p \geq (m+p)!$.

In \cref{sec:K2mresults}, we turn our attention to $K_{2*m}$. We begin by introducing some notation and developing tools to make progress toward showing $K_{2*m}$ is weakly enumeratively chromatic-choosable for each $m \in \N$.  Among these tools is an important recursive relation that is applied in \cref{sec:examples}.
\begin{lemma}\label{lem:lboundlistcolorfunction}
For any $m \geq 2$, \[ P_\ell \big(K_{2*m}, m\big)  \geq \left\lceil\left(\frac{m}{m-1} \right)^{m-1} P_\ell \big(K_{2*(m-1)}, m-1\big)\right\rceil.\]
\end{lemma} 
\noindent We conclude \cref{sec:K2mresults} by using these tools to prove that $\K{3}$ and $\K{4}$ are weakly enumeratively chromatic-choosable in \cref{prop:K222wecc,prop:K2222wecc}, respectively.

Finally, in \cref{sec:examples}, we give some nontrivial applications of \cref{thm:joincomplete}.  Note that since $K_{2*4}$ is weakly enumeratively chromatic-choosable, \cref{lem:lboundlistcolorfunction} tells us $P_{\ell}(K_{2*5},5) \geq 59$. We use this fact and \cref{cor:joincomplete} to show $K_p \vee K_{2*5}$ is weakly enumeratively chromatic-choosable whenever $p \geq 15$.  We give a similar result for $K_{2,2,3}$ after proving $P_{\ell}(K_{2,2,3},3) \geq 4$.

\section{Some Initial Tools and Joining Vertices}\label{sec:prelimresults}
In this section, we first introduce some tools that will be useful in subsequent sections. We conclude the section with a proof of \cref{thm:joincomplete}.
In \cite[Lemma 25]{AM25}, Allred and Mudrock proved the following to aid in proving graphs are enumeratively chromatic-choosable. This turns out to be rather useful in many of the proofs below.
\begin{lemma} \label{lem:prob} \cite[Lemma 25]{AM25}
Suppose $G$ is a graph with $V(G) = \{v_1, \ldots, v_n \}$.  Suppose that $m \in \N$, $d_1, \ldots, d_n$ are nonnegative integers, and $L$ is a list assignment for $G$ such that $|L(v_i)|=m+d_i$.  Then,
\[P(G,L) \geq \left \lceil \frac{P_{\ell}(G,m)\prod_{i=1}^n (m+d_i)}{m^n} \right \rceil.\]
\end{lemma}

In some situations, we use the following idea from matching theory to count the number of proper $L$-colorings.

\begin{lemma}\label{lem:matching} \cite{KW22}
Suppose $G$ is a bipartite graph with bipartition $X$, $Y$. Let $n = |X|$, and suppose the degree of each vertex in $X$ is at least $k$.  If $|N(S)| \geq |S|$ for every $S \subseteq X$, then the number of matchings in $G$ that saturate $X$ is at least $\prod_{i=1}^{\min \{k,n\}} (k-i+1)$. 
\end{lemma}
In order to adapt \cref{lem:matching} to counting proper $L$-colorings, we need to define a new bipartite graph based on our initial graph and the particular list assignment.
Let $L$ be an arbitrary $m$-assignment of a graph $G$.
The \emph{$L$-matching graph of $G$}, denoted by $G_{L}$, is the bipartite graph with partite sets $X=V(G)$, $Y=\{y~:~y\in L(v) \text{ for some }v\in V(G)\}$, and $G_{L}$ has an edge from a vertex $x\in X$ to a vertex $y\in Y$ if and only if $y\in L(x)$. A matching in $G_{L}$ saturating $X$ naturally corresponds to a proper $L$-coloring of $G$.

\begin{cor}\label{cor:matching}
    Suppose $G$ is an $n$-vertex graph and $L$ is a $k$-assignment for $G$ with $ k \leq n$.  If $| \bigcup_{v \in S} L(v) | \geq |S|$ for every $S \subseteq V(G)$, then $P(G,L) \geq k!$.
\end{cor}

\begin{proof}
    Suppose that $\M$ is the set of all matchings in $G_{L}$ that saturate $V(G)$. Also, assume that $\mathcal{L}$ is the set of all proper $L$-colorings of $G$. Consider the function $f : \M \rightarrow \mathcal{L}$ given by $f(M) = h$ where $h$ is the proper $L$-coloring that maps each $v \in V (G)$ to its adjacent vertex in $M$. Clearly, $f$ is injective, which means $|\M| \le |\mathcal{L}| = P (G, L)$.
    
    Note that $G_L$ is bipartite, $\deg_{G_L}(v)=k$ for each $v\in V(G)$, and $k\le n$.  
    Since $| \bigcup_{v \in S} L(v) | \geq |S|$ for every $S \subseteq V(G)$, \cref{lem:matching} implies that \[|\M|\ge \prod_{i=1}^{\min\{n,k\}} (k-i+1)=\prod_{i=1}^{k} (k-i+1)=k!.\]

    Hence, $k!\le |\M|\le P(G,L)$, as desired.
\end{proof}

We are now ready to prove \cref{thm:joincomplete}; however, we need a notion that will be used for the remainder of the paper.  When $G$ is a graph and $L$ is a list assignment for $G$, for each $c \in \bigcup_{v \in V(G)} L(v)$, the \emph{multiplicity of $c$ in $L$}, denoted by $\mu_L^G(c)$, is $|\{v \in V(G) : c \in L(v) \}|$.  When the list assignment and graph are clear from context, we simply write $\mu(c)$ for $\mu_L^G(c)$. 

\thmjoincomplete*
\begin{proof}
Our proof is by induction on $p$.  When $p=0$, we use the convention that $K_0 \vee G = G$. So, the desired result clearly holds when $p=0$.  Now, suppose $p>0$ and the desired result holds for all nonnegative integers less than $p$.

Let $H= K_p \vee G$, and suppose $\{u_1, \dots, u_p\}$ is the vertex set of the copy of $K_p$ used to form $H$. 
Let $A_{m+i} = \{u_i\}$ for each $i \in [p]$. Suppose $L$ is an $(m+p)$-assignment of $H$ satisfying $P(H,L)=P_{\ell}(H,m+p)$. 

Suppose $L(u_p)=\{c_1,\dots,c_{m+p}\}$. Without loss of generality, assume $\mu (c_1)\geq \mu (c_2)\geq \dots \geq \mu (c_{m+p})$. Let $n= |V(G)|$. We first prove that if $\mu(c_i)=n+p$ for each $i \in [p+1]$, $P(H,L)\geq (m+p)!$.

Let $\mathcal{P}$ be the set of $(p+1)$-tuples that consist of distinct elements of $[m+p]$. If $\textbf{x}=(x_1,\dots,x_{p+1}) \in \mathcal{P}$, then let $\mathcal{C}_{\textbf{x}}$ be the set of proper $L$-colorings of $H$ that use $c_i$ to color all the vertices of $A_{x_i}$ for each $i \in [p+1]$.  Let $H_{\textbf{x}}=H - \bigcup_{i=1}^{p+1} A_{x_i} $ and let $L_{\textbf{x}}$ be the $(m-1)$-assignment for $H_{\textbf{x}}$ given by $L_{\textbf{x}}(v) = L(v) - \{c_1, \dots, c_{p+1}\}$ for each $v \in V(H_{\textbf{x}})$.  Notice $|\mathcal{C}_{\textbf{x}}|=P(H_{\textbf{x}},L_{\textbf{x}})$.  By the hypotheses of the theorem and \zrefWithTheorem{examplesa,examplese}, $H_{\textbf{x}}$ is weakly enumeratively chromatic-choosable. Since $H_{\textbf{x}}$ is also uniquely $(m-1)$-colorable, $P(H_{\textbf{x}},L_{\textbf{x}}) \geq (m-1)!$.

Notice that the sets $\mathcal{C}_{\textbf{x}}$ are pairwise disjoint as $\textbf{x}$ ranges over $\mathcal{P}$, since any two distinct color classes $A_i$ and $A_j$ contain adjacent vertices. Since $|\mathcal{P}| = (m+p)!/(m-1)!$,
\begin{align*}
P(H,L) \geq \sum_{\textbf{y} \in \mathcal{P}} |\mathcal{C}_{\textbf{y}} |= \sum_{\textbf{y} \in \mathcal{P}} P(H_{\textbf{y}}, L_{\textbf{y}}) \geq |\mathcal{P}|(m-1)! = (m+p)!.
\end{align*}

Now, we assume $\mu(c_{p+1}) < n+p$.  For each $i \in [m+p]$, let $\mathcal{C}_i$ be the set of proper $L$-colorings of $H$ that color $u_p$ with $c_i$.  Clearly, $P(H,L) = \sum_{i=1}^{m+p} |\mathcal{C}_i|$.  Let $H'=H-\{u_p\}$ and let $L_i$ be the list assignment for $H'$ defined by $L_i(v)=L(v)-\{c_i\}$ for each $v\in V(H')$. Note that $H'=K_{p-1}\vee G$ and $|\mathcal{C}_i|=P(H',L_i)$.

For $i \in [p]$, the induction hypothesis implies $|\mathcal{C}_i| = P(H', L_i) \ge \min\{a_{p-1}, (m+p-1)!\}$. Notice for each $i \in [p+1:m + p]$, $L_i$ has at least one output of size $m + p$ since $\mu (c_i) < n + p$. So by \cref{lem:prob} and the induction hypothesis, $|\mathcal{C}_i| = P(H', L_i) \ge \min\{a_{p-1}, (m+p-1)!\} \left( \frac{m+p}{m+p-1} \right)$ for each $i \in [p+1:m + p]$.  This means
\[\sum_{i=1}^{m+p} |\mathcal{C}_i| \geq p \min\{a_{p-1}, (m+p-1)!\} + m \min\{a_{p-1}, (m+p-1)!\} \left( \frac{m+p}{m+p-1} \right)\]
Note that when $\min\{a_{p-1}, (m+p-1)!\}=(m+p-1)!$,  $\sum_{i=1}^{m+p} |\mathcal{C}_i| \geq (m+p)!$.  When $\min\{a_{p-1}, (m+p-1)!\}=a_{p-1}$, $\sum_{i=1}^{m+p} |\mathcal{C}_i| \geq \left( m \left( \frac{m+p}{m+p-1}\right) +p \right)a_{p-1}=a_p$.  This completes the induction step.
\end{proof}

\begin{cor} \label{cor:joincomplete}
Suppose $G$ is a uniquely $m$-colorable graph where $m \geq 2$ and the unique partition of $V(G)$ into $m$ independent sets is $\{A_i : i \in [m]\}$.  
Suppose $G$ is such that $G - \bigcup_{i=1}^{m-1} A_{x_i}$ is weakly enumeratively chromatic-choosable whenever $(x_1, \ldots, x_{m-1}) \in [m]^{m-1}$ and $P_{\ell}(G,m) \geq a_0$.  If there is an $N \in \N$ such that

\[\frac{a_0}{m!}\prod_{p=1}^N \left(1+ \frac{m}{(m-1+p)(m+p)}\right) \geq 1,\]
then 
$K_q \vee G$ is weakly enumeratively chromatic-choosable whenever $q \geq N$.
\end{cor}
\begin{proof}
Note that $G$ satisfies the hypotheses of \cref{thm:joincomplete}.  Consequently, using the notation of \cref{thm:joincomplete} for each $p \in \N$, we know that $P_{\ell}(K_p \vee G, m+p) \geq \min \{a_p, (m+p)! \}$.

If there is a $d < N$ such that $a_d \geq (m+d)!$, we have that $K_d \vee G$ is weakly enumeratively chromatic-choosable.  Then, \zrefWithTheorem{examplese} implies $K_q \vee G$ is weakly enumeratively chromatic-choosable whenever $q \geq d$.

So, we may assume $a_p < (m+p)!$ for all $p \in [N-1]$.  By the recursive formula in \cref{thm:joincomplete}, $a_N =  \prod_{p=1}^{N}\left( m \left( \frac{m+p}{m+p-1}\right) +p \right)a_0$. Thus by our assumption on $a_0$ above we obtain 
\begin{align*}
a_N &\geq  m! \prod_{p=1}^{N} (m+p) \left(  \left( \frac{m}{m+p-1}\right) +\frac{p}{m+p} \right)\left( \frac{(m + p -1)(m+p)}{(m + p - 1)(m + p) +m } \right) \\
&= m! \prod_{p=1}^{N} (m+p)\left( \frac{m(m+p)+p(m+p-1)}{(m+p-1)(m+p)} \right)\left( \frac{(m + p -1)(m+p)}{(m + p - 1)(m + p) +m } \right) \\
&=(m+N)!.
\end{align*}
So, \cref{thm:joincomplete} implies $K_N \vee G$ is weakly enumeratively chromatic-choosable. Hence, \zrefWithTheorem{examplese} implies the desired result.
\end{proof}

\section{Results for $K_{2*m}$}\label{sec:K2mresults}
We begin by establishing some notational conventions that will be used throughout this section. Whenever $G=\K{m}$ for some $m \geq 2$, we will assume the partite sets of $G$ are $A_1,\dots, A_m$ where $A_i=\{x_i,y_i\}$ for each $i\in [m]$. Additionally, when $L$ is a list assignment for such a $G$ we let $\C{i,c}$ be the set of proper $L$-colorings of $G$ that color $x_i$ and $y_i$ with $c$ (note such a set may be empty).

Before we show that $\K{m}$ is weakly enumeratively chromatic-choosable when $m \in \{3,4\}$, we prove three lemmas that we will need.

\begin{lemma}\label{lem:multhc1}
Suppose $m \in \N$ and $G=K_{2*m}$.  Let $L$ be an $m$-assignment of $G$. Suppose that $|L(x_i)\cap L(y_i)|\le m-i$ for all $i\in [m]$. Then, $P(G,L)\ge m!$.
\end{lemma}

\begin{proof}
    The strategy for the proof is to use \cref{cor:matching} which says that the desired result will follow if $\abs{\bigcup_{v \in S} L(v)} \geq |S|$ whenever $S \subseteq V(G)$.  Suppose $S \subseteq V(G)$.  If $S = \emptyset$, we clearly have $|\bigcup_{v \in S} L(v)| =0= |S|$.  Also, if $S$ is nonempty and $|S| \leq m$, we have $\abs{\bigcup_{v \in S} L(v)} \geq m \geq |S|$ since $L$ is an $m$-assignment.  So, we may suppose that $|S| > m$.

    By the Pigeonhole Principle, there is a $t \in [m]$ such that $A_t \subseteq S$.  Let $q$ be the largest element in $[m]$ such that $A_q \subseteq S$.  By the definition of $q$, we note that $|S| \leq m+q$.  Let $d = |S|-m$ so that $|S|=m+d$, and notice $d \leq q$.  
    
    Since $|L(x_q) \cap L(y_q)| \leq m-q$ and $L$ is an $m$-assignment, we have that $|L(x_q) \cup L(y_q)| \geq m+q$.  Consequently,
    \[\left|\bigcup_{v \in S} L(v)\right| \geq |L(x_q) \cup L(y_q)| \geq m+q \geq m+d = |S|\]
    which completes the proof.
\end{proof}

\begin{customlem}{\ref{lem:lboundlistcolorfunction}}
\label{lem:lboundlistcolorfunctionintro}
Suppose $m \geq 2$, $G= K_{2*m}$, and $L$ is an $m$-assignment for $G$.
Let $z\in \bigcup_{v\in V(G)} L(v)$ be such that
$\mu(z)\geq \mu(c)$ for each
$c\in \bigcup_{v\in V(G)} L(v)$.
Then the following statements hold.
\begin{statement}
    \item\label{bound1intro}
    If $\mu(z)\leq m$, then $P(G,L)\geq m!$.
    
    \item\label{bound2intro}
    If $\mu(z)=m+q$ for $q\in [m]$, then
    \[
    P(G,L)\geq
    q \left \lceil \left(\frac{m}{m-1}\right)^{m-q}
    P_{\ell}\big(\K{(m-1)},m-1\big) \right \rceil.
    \]
\end{statement}
Consequently,
\[
P_\ell\big(K_{2*m},m\big)
\geq
\left\lceil
\left(\frac{m}{m-1}\right)^{m-1}
P_\ell\big(K_{2*(m-1)},m-1\big)
\right\rceil.
\]
\end{customlem}

\begin{proof}
First, we prove \ref{bound1intro}.  Suppose $\mu(z)\le m$ for each $z\in \bigcup_{v\in V(G)} L(v)$.
Consider the $L$-matching graph of $G$, $G_L$, with bipartition $X=V(G)$, $Y$.  We will show the desired result by using \cref{cor:matching}.
Suppose $S\subseteq V(G)$. Let  $A$ be the set of edges in $G_L$ with an endpoint in $S$, and let $B$ be the set of edges in $G_L$ with an endpoint in $N(S)$. 
Since each vertex in $X$ has degree $m$, we have that $|A|=m|S|$. 
Now, since $\mu (z)\le m$, it follows that the degree of each vertex in $Y$ is at most $m$. 
Thus, $|B| \leq m|N(S)|$. Since $A\subseteq B$, it follows that $m|S| = |A| \leq |B| \leq m|N(S)|$ which means $|S| \leq |N(S)|$.  
Finally, since $|N(S)| = |\bigcup_{v \in S} L(v)|$, the hypothesis of \cref{cor:matching} is satisfied.  Therefore, $P(G,L)\ge m!$ as desired.

We will use \cref{lem:prob} to prove \ref{bound2intro}. 
Since $\mu(z)=m+q$, without loss of generality assume that $z\in L(x_i)$ and $z\in L(y_i)$ for $i\in [q]$.  Note that $\C{1,z}, \dots, \C{q,z}$ are pairwise disjoint.
    So $P(G,L)\ge \sum\limits_{i=1}^{q} |\C{i,z}|$.
    We will now establish a lower bound on $|\C{i,z}|$ for each $i \in [q]$. Suppose for some $j\in [q]$ that $x_j$ and $y_j$ are both colored with $z$. Then, let $G'=G-A_j$ and $L'(v)=L(v)-\{z\}$ for each $v \in V(G')$. Note that $|\C{j,z}|=P(G',L')$ and $|L'(v)|\ge m-1$ for each $v\in V(G')$. Moreover, exactly $m-q$ of the lists associated with $L'$ have size $m$.
    
    \cref{lem:prob} implies that
    \[|\C{j,z}|=P(G',L')\ge \left\lceil \frac{P_{\ell}(G',m-1)m^{m-q} (m-1)^{m+q-2}}{(m-1)^{2m-2}}\right\rceil=\left\lceil \frac{P_{\ell}(G',m-1)m^{m-q}}{(m-1)^{m-q}}\right\rceil.\]

    Hence, 
    \[P(G,L)\ge\sum\limits_{i=1}^{q} |\C{i,z}|\ge q\left\lceil \frac{P_{\ell}(G',m-1)m^{m-q}}{(m-1)^{m-q}}\right\rceil\ge  q \left \lceil \left(\frac{m}{m-1}\right)^{m-q} P_{\ell}\big(\K{(m-1)},m-1\big) \right \rceil \] as desired.
    
For the final inequality, suppose $t\in[m-1]$.  Then, $m-t \geq 1$ which implies $mt+m-t-1 \geq mt$ which implies 
\[\frac{m-1}{m} \geq \frac{t}{t+1}
\hspace{.25cm}\text{ which implies }\hspace{.25cm}
(t+1)\left(\frac{m}{m-1}\right)^{m-t-1}  
 \geq t\left(\frac{m}{m-1}\right)^{m-t}.\]

Since this holds for each $t\in [m-1]$, it follows that for each $q\in[m]$
\[q\left(\frac{m}{m-1}\right)^{m-q} \geq \left(\frac{m}{m-1}\right)^{m-1}.\]

Also, note that in the case $\mu(z) \leq m$, \[P(G,L)\geq m!=m(m-1)!=mP\big(\K{(m-1)},m-1\big)\ge \left(\frac{m}{m-1}\right)^{m-1}P_{\ell}\big(\K{(m-1)},m-1\big).\]  These two facts imply the final inequality.
\end{proof}

\begin{lemma}\label{lem:mult2m-t}
Suppose $t \in \{2,3\}$, $m\ge t+1$, $G=K_{2*m}$, and $L$ is an $m$-assignment of $G$. Suppose there exists a $z\in \bigcup_{v\in V(G)} L(v)$ such that $\mu(z)=2m-t$ and $\mu(c)\leq 2m-t$ for each $c \in \bigcup_{v \in V(G)} L(v)$.
Suppose $z\in L(x_i)\cap L(y_i)$ for each $i\in [m-t]$.
Then $P(G,L)\ge m P_{\ell}(\K{(m-1)},m-1)$ if at least one of the following statements hold.
\begin{statement}
    \item\label{statement1} $z \in L(x_s) \cap L(y_s)$ for some $s \in [m-t+1 : m]$.
    \item\label{statement2} $|L(x_s)\cap L(y_s)| \geq m-t+1$ for some $s\in [m-t]$.
    \item\label{statement3} $L(x_s)\cap L(y_s) \neq \emptyset$ for some $s\in [m-t+1 : m]$.
\end{statement}
\end{lemma}

\begin{proof}
The strategy for the proof is as follows.  First, we show the desired result holds when \ref{statement1} holds.  Second we show the desired result holds when \ref{statement1} doesn't hold and \ref{statement2} does hold.  Finally, we show the desired result holds when neither \ref{statement1} nor \ref{statement2} holds, but \ref{statement3} does hold.

First, suppose that $z\in L(x_s)\cap L(y_s)$ for some $s\in [m-t+1:m]$.  Without loss of generality, we may assume that $z\in L(x_{m-t+1})\cap L(y_{m-t+1})$.
Since $\C{1,z}, \dots, \C{m-t+1,z}$ are pairwise disjoint, it follows that $P(G,L)\ge\sum\limits_{i=1}^{m-t+1} |\C{i,z}|$. 

Now, suppose $q$ is a fixed element of $[m-t+1]$.  We wish to establish a lower bound on $|\C{q,z}|$ by using \cref{lem:prob}.  Suppose that $x_q$ and $y_q$ are both colored with $z$.
Let $G'=G-A_q$, and note $G' = K_{2*(m-1)}$.  Also, let $L'(v)=L(v)-\{z\}$ for each $v \in V(G')$, and note $|\C{q,z}|=P(G',L')$. Since $\mu(z)=2m-t$, $L'$ maps exactly $t$ elements of $V(G')$ to a set of size $m$, and $L'$ maps the remaining elements of $V(G')$ to a set of size $m-1$.  So, by \cref{lem:prob}, we have that $P(G',L')\ge P_{\ell}(K_{2*(m-1)},m-1)(1+\frac{1}{m-1})^t$. It follows that
\begin{equation}
    \begin{split}
    \sum\limits_{i=1}^{m-t+1}|\C{i,z}|&\ge (m-t+1)\left(1+\frac{1}{m-1}\right)^t P_{\ell}(K_{2*(m-1)},m-1)\\
    &\geq (m-t+1)\left(1+\frac{t}{m-1}\right)P_{\ell}(K_{2*(m-1)},m-1) \\
    &= \left(m-t+1 + \frac{t(m-1)+t(-t+2)}{m-1}\right)P_{\ell}(K_{2*(m-1)},m-1)\\
    &\geq \left(m + 1 + \frac{t(-t+2)}{t} \right)P_{\ell}(K_{2*(m-1)},m-1) \\
    &\ge m P_{\ell}(K_{2*(m-1)},m-1). \\
    \end{split}
    \label{calc:3}
\end{equation}

 We now turn our attention to showing the desired result holds when \ref{statement1} doesn't hold and \ref{statement2} does hold (note that we don't need to use the fact \ref{statement1} doesn't hold for this part of the proof).  Suppose that $|L(x_s)\cap L(y_s)|\ge m-(t-1)$ for some $s\in [m-t]$.  Suppose $L(x_s)\cap L(y_s)\supseteq \{c_1,\dots, c_{m-t+1}\}$.
Since $\C{s, c_1}, \C{s, c_2}, \dots , \C{s, c_{m-t+1}}$ are pairwise disjoint, it follows that $P(G,L)\ge\sum\limits_{j=1}^{m-t+1} |\C{s,c_j}|$.  Now, suppose $q$ is a fixed element of $[m-t+1]$.  We will establish a lower bound on $|\C{s,c_q}|$ by using \cref{lem:prob}.  Suppose that $x_s$ and $y_s$ are both colored with $c_q$. Let $G'=G-A_s$, and note $G' = K_{2*(m-1)}$.  Also, let $L'(v)=L(v)-\{c_q\}$ for each $v \in V(G')$, and note $|\C{s,c_q}|=P(G',L')$.  Since $\mu(c_q) \leq 2m-t$, $L'$ maps at least $t$ elements of $V(G')$ to a set of size $m$ and each element of $V(G')$ to a set of size at least $m-1$. So, by \cref{lem:prob}, we have that $P(G',L')\ge 
P_{\ell}(K_{2*(m-1)},m-1)(1+\frac{1}{m-1})^t$.  Using the ideas in~(\ref{calc:3}) above, we get \[\sum\limits_{j=1}^{m-t+1} |\C{s,c_j}|\ge (m-t+1)\left(1+\frac{1}{m-1}\right)^t P_{\ell}(K_{2*(m-1)},m-1)\ge m P_{\ell}(K_{2*(m-1)},m-1)\]  
as desired. 

Finally, suppose that \ref{statement1} and \ref{statement2} do not hold and $L(x_s)\cap L(y_s)\ne \emptyset$ for some $s\in [m-t+1 : m]$.  Suppose $c\in L(x_s)\cap L(y_s)$.  Without loss of generality, we may assume that $z\in L(x_{i})$ and $z \notin L(y_i)$ for each $i \in [m-t+1:m]$.  

For each $i \in [m-t]$ let $\C{i,z}'$ be the set of proper $L$-colorings of $G$ that color $x_i$ and $y_i$ with $z$, but does not color $y_s$ with $c$.  Note that
$\C{1,z}',\dots,\C{m-t,z}',\C{s,c}$ are pairwise disjoint. 
We will now establish a lower bound on $|\C{q,z}'|$ when $q$ is a fixed element of $[m-t]$ by using \cref{lem:prob}.  Suppose $x_q$ and $y_q$ are both colored with $z$. Let $G'=G-A_q$, and note $G'=\K{(m-1)}$.  Also, let $L'(v)=L(v)-\{z\}$ for each $v \in V(G') - \{y_{s}\}$ and $L'(y_{s}) = L(y_{s}) - \{c\}$, and note $|\C{q,z}'|=P(G',L')$.  By the definition $L'$, we have that $L'$ maps the vertices in $(\{y_j:j \in [m-t+1:m]\} - \{y_s\})$ to a set of size $m$, and $L'$ maps all other vertices of $G'$ to a set of size $(m-1)$.  So, by \cref{lem:prob}, we have that $|\C{q,z}'|\ge
P_{\ell}(K_{2*(m-1)},m-1)(1+\frac{1}{m-1})^{t-1}$.

Now, we will establish a lower bound on $|\C{s,c}|$. Using the same argument as above for the lower bound on $|\C{s,c_q}|$, we get
\[
|\C{s,c}| \ge P_\ell(K_{2\ast(m-1)},m-1)\left(1+\frac{1}{m-1}\right)^t.
\]
Putting all this together, we obtain
\begin{align*}
&|\C{s,c}|+\sum_{i=1}^{m-t}|\C{i,z}'|\\
&\ge \left(\frac{P_\ell(K_{2\ast(m-1)},m-1)m^t}{(m-1)^t}+(m-t)\frac{P_\ell(K_{2\ast(m-1)},m-1)m^{t-1}}{(m-1)^{t-1}}\right)\\
&\ge \left(\frac{m^t}{(m-1)^t}+\frac{(m-1)(m-t)m^{t-1}}{(m-1)^t}\right)P_\ell(K_{2\ast(m-1)},m-1)\\
&\ge m\left(\frac{m^{t+1}-tm^t+tm^{t-1}}{m(m-1)^t}\right)P_\ell(K_{2\ast(m-1)},m-1).
\end{align*}
It is easy to check that $m^{t+1}-tm^t+tm^{t-1}>m(m-1)^t$ when $t\in\{2,3\}$. Consequently,
\[
P(G,L)\ge |\C{s,c}|+\sum_{i=1}^{m-t}|\C{i,z}'|\ge mP_\ell(K_{2\ast(m-1)},m-1)
\]
as desired.
\end{proof}

For the rest of the section, when $G=\K{m}$ and $L$ is a list assignment for $G$, we let $d_i=|L(x_i)\cap L(y_i)|$ for each $i\in[m]$. Whenever we apply \cref{lem:multhc1}, we will assume, without loss of generality, that $d_1\geq d_2\geq\cdots\geq d_m$.

\begin{pro}\label{prop:K222wecc}
$K_{2*3}$ is weakly enumeratively chromatic-choosable. 
\end{pro}
\begin{proof}
     Suppose $G=K_{2*3}$. We know that $P_{\ell}(G,3) \leq P(G,3) = 6$.  Suppose $L$ is an arbitrary $3$-assignment for $G$. We must show that $6\leq P(G,L)$.  We will consider cases based on the maximum multiplicity of colors in $\bigcup_{v\in V(G)} L(v)$. Suppose that $z$ is such a color of maximum multiplicity.
     
Suppose $\mu(z)=6$ for some $z\in \bigcup_{v\in V(G)}L(v)$. Then, using \zrefWithTheorem{bound2intro} and the fact that $K_{2*2}$ is enumeratively chromatic-choosable by \zrefWithTheorem{examplesb}, we have that $P(G,L)\ge 3 P_{\ell}(K_{2*2},2)=(3)(2)=6$; the conclusion follows.
        
Next, suppose $\mu(z)=5$. Then \zrefWithTheorem{bound2intro} and \zrefWithTheorem{examplesb} imply that $P(G,L)\ge 6$;  the conclusion follows.

Next, assume that $\mu(z)=4$. If at least one of the three statements in \cref{lem:mult2m-t} hold, then $P(G,L)\ge 3 P_{\ell}(\K{2},2)=6$ by \cref{lem:mult2m-t} and \zrefWithTheorem{examplesb}. So, we may assume the three statements in \cref{lem:mult2m-t} do not hold; specifically, we have $d_1 \leq 1$ and $d_2=d_3=0$.   So, $d_1 \leq 1 \leq 2$, $d_2=0 \leq 1$, and $d_3=0$.  Thus, \cref{lem:multhc1} implies the desired.

     We may therefore assume that $\mu(c)\le 3$ for all $c\in \bigcup_{v\in V(G)}L(v)$.  Then \zrefWithTheorem{bound1intro} implies that $P(G,L)\ge 3!=6$;  the conclusion follows.
\end{proof}

\begin{pro}\label{prop:K2222wecc}
    $\K{4}$ is weakly enumeratively chromatic-choosable.
\end{pro}

\begin{proof}
Suppose $G=K_{2*4}$. We know that $P_{\ell}(G,4)\leq P(G,4)=24$. 
    Suppose $L$ is an arbitrary 4-assignment for $G$. We need to show that $24 \leq P(G,L)$.  We will consider cases based on the maximum multiplicity of colors in $\bigcup_{v\in V(G)} L(v)$. Suppose that $z$ is such a color of maximum multiplicity. 
    
    Suppose $\mu(z)=8$.  Then, using \zrefWithTheorem{bound2intro} and the fact that $K_{2*3}$ is weakly enumeratively chromatic-choosable by \cref{prop:K222wecc}, we have that $P(G,L)\ge 4 P_{\ell}(\K{3},3)=(4)(6)=24$; the conclusion follows.
    
    Next, suppose $\mu(z)=7$. Then \zrefWithTheorem{bound2intro,prop:K222wecc} imply that $P(G,L)\ge 24$;  the conclusion follows.

    Next, assume that $\mu(z)=6$.  If at least one of the three statements in \cref{lem:mult2m-t} hold, then $P(G,L)\ge 4 P_{\ell}(\K{3},3)=24$ by \cref{lem:mult2m-t,prop:K222wecc}. So, we may assume the three statements in \cref{lem:mult2m-t} do not hold; specifically, we have $d_1 \leq 2$, $d_2 \leq 2$, and $d_3=d_4=0$.   So, $d_1 \leq 2 \leq 3$, $d_2\leq 2 \leq 2$, $d_3=0\leq 1$, and $d_4=0$. \cref{lem:multhc1} implies the desired conclusion. 

    Next, assume that $\mu(z)=5$. If at least one of the three statements in \cref{lem:mult2m-t} hold, then $P(G,L)\ge 4 P_{\ell}(\K{3},3)=24$ by \cref{lem:mult2m-t,prop:K222wecc}. So, we may assume the three statements in \cref{lem:mult2m-t} do not hold; specifically, we have $d_1 \leq 1$ and $d_2=d_3=d_4=0$.   So, $d_1 \leq 1 \leq 3$, $d_2=0 \leq 2$, $d_3=0\leq 1$, and $d_4=0$. \cref{lem:multhc1} implies the desired conclusion.

    We may therefore assume that $\mu(c)\le 4$ for all $c\in \bigcup_{v\in V(G)}L(v)$.  Then \zrefWithTheorem{bound1intro} implies that $P(G,L)\ge 4!=24$;  the conclusion follows.
\end{proof}

In \cref{tab:lowerbounds} we apply \cref{lem:lboundlistcolorfunction} and use \cref{prop:K2222wecc} to get a lower bound on $P_{\ell}(\K{m},m)$ for each $m \in [5:15]$. One can compare the lower bounds to $m!$ which is the lower bound required to prove $\K{m}$ is weakly enumeratively chromatic-choosable.
\begin{table}[h]
    \centering
    \begin{tabular}{c|r|r}
        $m$ & Lower Bound & $m!$ \\
        \hline
        5  & 59      & 120 \\
        6  & 147     & 720 \\
        7  & 371     & 5040 \\
        8  & 945     & 40320 \\
        9  & 2425    & 362880 \\
        10 & 6260    & 3628800 \\
        11 & 16237   & 39916800 \\
        12 & 42285   & 479001600 \\
        13 & 110493  & 6227020800 \\
        14 & 289559  & 87178291200 \\
        15 & 760716  & 1307674368000 \\
    \end{tabular}
    \caption{Lower bounds of $P_{\ell}(\K{m}, m)$ compared to $m!$}
    \label{tab:lowerbounds}
\end{table}
\FloatBarrier

\section{Applying \cref{thm:joincomplete}}\label{sec:examples}
In this section, we apply \cref{thm:joincomplete} to $\K{5}$ and $K_{2,2,3}$ via \cref{cor:joincomplete}.  We begin with $\K{5}$.

\begin{example}\label{ex:k25}
Let $G = \K{5}$. We know that $G$ is uniquely 5-colorable. Suppose the unique partition of $V(G)$ into $5$ independent sets is $\{A_i:i\in [5]\}$. 
Furthermore, by \zrefWithTheorem{examplesa,examplesb} and \cref{prop:K222wecc,prop:K2222wecc}, $G - \bigcup_{i=1}^{4} A_{x_i}$ is weakly enumeratively chromatic-choosable for all $(x_1, x_2, x_3, x_4) \in [5]^4$. By \cref{lem:lboundlistcolorfunction,prop:K2222wecc}, $P_\ell \big(K_{2*5}, 5\big) \geq \left(\frac{5}{4} \right)^{4} P_\ell \big(K_{2*(4)}, 4\big) = \left(\frac{5}{4}\right)^4 \cdot 24$.  Consequently, $P_{\ell}(K_{2*5},5) \geq 59$.

Note that \[\frac{59}{5!}\prod_{p=1}^{15} \left(1+ \frac{5}{(4+p)(5+p)}\right) \geq 1.\]  Consequently, \cref{cor:joincomplete} implies that $K_m \vee K_{2*5}$ is weakly enumeratively chromatic-choosable whenever $m \geq 15$.
\end{example}

To apply \cref{cor:joincomplete} to $K_{2,2,3}$ we first need an $a_0$ such that $P_{\ell}(K_{2,2,3},3) \geq a_0$.

\begin{pro}\label{prop:K223wecc}
$P_{\ell}(K_{2,2,3},3)\geq 4$. 
\end{pro}
\begin{proof}
    Suppose $G=K_{2,2,3}$. Suppose $L$ is an arbitrary $3$-assignment for $G$. We wish to show $P_{\ell}(K_{2,2,3},3)\geq 4$. Suppose the tripartition is $A_1=\{ x_1,y_1,z_1\}$, $A_2=\{x_2,y_2\}$, and $A_3=\{x_3,y_3\}$. Let $\C{i,q}$ be the set of proper $L$-colorings of $G$ that color each vertex of $A_i$ with $q$. Note that $\C{1,q}, \C{2,q}, \C{3,q}$ are pairwise disjoint for any color $q$.  We will consider cases based on the maximum multiplicity of colors in $\bigcup_{v\in V(G)} L(v)$.  Suppose $z$ is such a color of maximum multiplicity.

    Suppose $\mu(z)=7$. Consider $\C{1,z}$ and let $G'=G-A_1$ and $L'(v)=L(v)-\{z\}$ for each $v\in V(G')$. Note that $G'$ is a copy of $C_4$. By \zrefWithTheorem{examplesb}, we know that $G'$ is enumeratively chromatic-choosable. Thus, $P_{\ell}(G',2)=2$ and $|\C{1,z}| \geq 2$. Consider $\C{2,z}$ and let $G'=G-A_2$ and $L'(v)=L(v)-\{z\}$ for each $v\in V(G')$. Note that $G'$ is a copy of $K_{2,3}$. By \zrefWithTheorem{examplesc}, we know that $G'$ is enumeratively chromatic-choosable. Thus, $P_{\ell}(G',2) = 2$ and $|\C{2,z}| \geq 2$. By a similar argument, we know that $|\C{3,z}|\geq 2$. Since $\C{1,z}, \C{2,z}, \C{3,z}$ are pairwise disjoint, $P(G,L)\geq 6$.

    Suppose $\mu(z)=6$. Let $A_i$ be such that $z \in L(v)$ for each $v \in A_i$. Note that there are two such $i\in[3]$ for which this is true. Let $G'=G-A_i$ and $L'(v)=L(v)-\{z\}$ for each $v\in V(G')$.  Notice that exactly one output of $L'$ is size 3 and the remaining outputs of $L'$ are size 2. Note that $G'$ is either a copy of $C_4$ or $K_{2,3}$ and therefore by \zrefWithTheorem{examplesb,examplesc}, $P_{\ell}(G',2)=2$. By \cref{lem:prob}, $|\C{i,z}|\geq (3/2)2 =3$. Because there are two such $i$, $P(G,L)\geq 6$.

Suppose $\mu(z)=5$. Let $A_i$ be such that $z\in L(v)$ for each $v\in A_i$. Note that there is at least one such $i\in[3]$ for which this is true. Let $G'=G-A_i$ and $L'(v)=L(v)-\{z\}$ for each $v\in V(G')$.  Notice that exactly two outputs of $L'$ are size 3 and the remaining outputs of $L'$ are size 2. Note that $G'$ is either a copy of $C_4$ or $K_{2,3}$ and therefore by \zrefWithTheorem{examplesb,examplesc}, $P_{\ell}(G',2)=2$. By \cref{lem:prob}, $|\C{i,z}|\geq \ceil{2(3/2)^2} =5$. Thus, $P(G,L)\geq 5$.

    Suppose $\mu(z)=4$. First, suppose there exist an $A_i$ and a color $q$ such that $q \in L(v)$ for each $v \in A_i$.  Let $G'=G-A_i$ and $L'(v)=L(v)-\{q\}$ for each $v\in V(G')$.  Notice that there are at least 3 outputs of $L'$ of size 3 and all the outputs of $L'$ are of size at least 2. Note that $G'$ is either a copy of $C_4$ or $K_{2,3}$ and therefore by \zrefWithTheorem{examplesb,examplesc}, $P_{\ell}(G',2)=2$. By \cref{lem:prob}, $|\C{i,q}|\geq \ceil{2(3/2)^3} =7$. Thus, $P(G,L)\geq 7$. So we may assume that $L(x_1)\cap L(y_1)\cap L(z_1)=\emptyset$, $L(x_2)\cap L(y_2)=\emptyset$, and $L(x_3)\cap L(y_3)=\emptyset$ and that without loss of generality $z$ is an element of $L(x_1), L(y_1), L(x_2)$, and $L(x_3)$. Color $x_1, y_1$ with $z$ and let $G''=G-\{x_1,y_1\}$ and $L''(v)=L(v)-\{z\}$ for each $v \in V(G'')$. Note that $G''$ is a copy of $K_1\vee C_4$.
    Notice that $P(G,L) \geq P(G'',L'')$.  We will now show that there are at least 4 proper $L''$-colorings of $G''$.  We split into cases based on the size of $L''(x_2) \cap L''(x_3)$. 

    Suppose $\left|L''(x_2) \cap L''(x_3)\right|=2$. Color $z_1$ with an element $q$ in $L''(z_1)-L''(x_2)$.  Let $L'''(v)=L''(v) - \{q\}$ for $v \in A_2 \cup A_3$.  Note $|L'''(v)| \geq 2$ for each such $v$.  Since $L'''(x_2)$ and $L'''(y_3)$ are disjoint along with $L'''(y_2)$ and $L'''(x_3)$, the number of proper $L''$-colorings of $G''$ that color $z_1$ with $q$ is at least
    \[\left|L'''(x_2)\right|\big(|L'''(x_3)|-1\big)\left|L'''(y_2)\right|\big(|L'''(y_3)|-1\big) \geq 2(1)(2)(1)=4.\]
    
    Suppose $\left|L''(x_2) \cap L''(x_3)\right|=1$. Without loss of generality, assume $L''(x_2) = \{a,b\}$ and $L''(x_3) = \{a,d\}$. Suppose $a\in L''(z_1)$ and color $z_1$ with $a$. Let $L'''(v)=L''(v) - \{a\}$ for $v \in A_2 \cup A_3$.  Then, we must color $x_2$ and $x_3$ with $b$ and $d$ respectively.  To complete a proper $L''$-coloring of $G''$, we must color $y_2$ with an element of $L'''(y_2) - \{d\}$ and then $y_3$ with an element of $L'''(y_3) - \{b\}$ that is different than the color used on $y_2$.  Thus, the number of proper $L''$-colorings of $G''$ that color $z_1$ with $a$ assuming $a\in L''(z_1)$ is at least 2. 
    
    Suppose $b\in L''(z_1)$ and color $z_1$ with $b$. Let $L'''(v)=L''(v) - \{b\}$ for $v \in A_2 \cup A_3$.  Note $|L'''(y_2)| = 3$ and $|L'''(y_3)| \geq 2$. Then, we must color $x_2$ and $x_3$ with $a$ and $d$ respectively. To complete a proper $L''$-coloring of $G''$, we must color $y_2$ with an element of $L'''(y_2) - \{d\}$ and then $y_3$ with an element in $L'''(y_3)$ that is different than the color used on $y_2$. Thus, the number of proper $L''$-colorings of $G''$ that color $z_1$ with $b$ assuming $b\in L''(z_1)$ is at least 2. By a similar argument, if $d\in L''(z_1)$, the number of proper $L''$-colorings of $G''$ that color $z_1$ with $d$ assuming $d\in L''(z_1)$ is at least 2.

    Suppose $q\in L''(z_1)$ and $q \notin \{a,b,d\}$. Color $z_1$ with $q$. Let $L'''(v)=L''(v) - \{q\}$ for $v \in A_2 \cup A_3$.  Notice that each output of $L'''$ is of size at least 2.  So, by \zrefWithTheorem{examplesb} there are at least 2 proper $L'''$-colorings of $G[A_2 \cup A_3]$. Thus, the number of proper $L''$-colorings of $G''$ that color $z_1$ with $q$ assuming $q\in L''(z_1)$ is at least 2.

    This shows that for any color $s\in L''(z_1)$, there are at least 2 proper $L''$-colorings of $G''$ that color $z_1$ with $s$.  Consequently, there are at least $6$ proper $L''$-colorings of $G''$ when $|L''(x_2) \cap L''(x_3)|=1$.

     Finally, suppose $\left|L''(x_2) \cap L''(x_3)\right|=0$. Assume $L''(x_2) = \{a,b\}$ and $L''(x_3) = \{d,r\}$. Suppose $a\in L''(z_1)$ and color $z_1$ with $a$. Let $L'''(v)=L''(v) - \{a\}$ for $v \in A_2 \cup A_3$.  Then, we must color $x_2$ with $b$. Then, color $y_3$ with an element $t$ of $L'''(y_3) - \{b\}$.  To complete a proper $L''$-coloring of $G''$, we must properly color the edge $y_2x_3$ so that the color used on $y_2$ is in $L'''(y_2)-\{t\}$ and the color used on $x_3$ is in $\{d,r\}$.  Since $|L'''(y_2)-\{t\}| \geq 2$, the number of proper $L''$-colorings of $G''$ that color $z_1$ with $a$ assuming $a\in L''(z_1)$ is at least 2. 

     In fact similar arguments can be used to show that for any $x \in  \{a,b,d,r\}$ if $x \in L''(z_1)$ the number of proper $L''$-colorings of $G''$ that color $z_1$ with $x$ is at least 2.

    Suppose $q\in L''(z_1)$ and $q \notin \{a,b,d,r\}$. Color $z_1$ with $q$. Let $L'''(v)=L''(v) - \{q\}$ for $v \in A_2 \cup A_3$.  Notice that each output of $L'''$ is of size at least 2.  So, by \zrefWithTheorem{examplesb} there are at least 2 proper $L'''$-colorings of $G[A_2 \cup A_3]$. Thus, the number of proper $L''$-colorings of $G''$ that color $z_1$ with $q$ assuming $q\in L''(z_1)$ is at least 2.

    This shows that for any color $s\in L''(z_1)$, there are at least 2 proper $L''$-colorings of $G''$ that color $z_1$ with $s$.  Consequently, there are at least $6$ proper $L''$-colorings of $G''$ when $|L''(x_2) \cap L''(x_3)|=0$.
    
Finally, assume that $\mu(c)\le 3$ for all $c\in \bigcup_{v\in V(G)} L(v)$.  Consider the $L$-matching graph of $G$, $G_L$, with bipartition $X=V(G)$, $Y$.  We will show that $P(G,L)\ge 6$ using \cref{cor:matching}.
Suppose $S\subseteq V(G)$. Let  $A$ be the set of edges in $G_L$ with an endpoint in $S$, and let $B$ be the set of edges in $G_L$ with an endpoint in $N(S)$. 
Since each vertex in $X$ has degree $3$, we have that $|A|=3|S|$. 
Now, since $\mu (c)\le 3$, it follows that the degree of each vertex in $Y$ is at most $3$ and that $|B| \leq 3|N(S)|$. Since $A\subseteq B$, it follows that $3|S| = |A| \leq |B| \leq 3|N(S)|$ which means that $|S| \leq |N(S)|$.  
Finally, since $|N(S)| = |\bigcup_{v \in S} L(v)|$, the hypothesis of \cref{cor:matching} is satisfied.  Therefore, $P(G,L)\ge 6$ as desired.
\end{proof}

Finally, we apply \cref{cor:joincomplete} to $K_{2,2,3}$.
\begin{example}\label{ex:k223}
Let $G = K_{2,2,3}$. We know that $G$ is uniquely 3-colorable. Suppose the tripartition is $A_1, A_2, A_3$. Furthermore, by \zrefWithTheorem{examplesa,examplesb,examplesc}, $G - \bigcup_{i=1}^{2} A_{x_i}$ is weakly enumeratively chromatic-choosable for all $(x_1, x_2) \in [3]^2$. By \cref{prop:K223wecc}, $P_{\ell}(K_{2,2,3},3)\geq 4$.  Note that \[\frac{4}{3!}\prod_{p=1}^{3} \left(1+ \frac{3}{(2+p)(3+p)}\right) \geq 1.\]  Consequently, \cref{cor:joincomplete} implies that $K_m \vee K_{2,2,3}$ is weakly enumeratively chromatic-choosable whenever $m \geq 3$.
\end{example}

\section*{Declaration of generative AI and AI-assisted technologies in the manuscript preparation process.} During the preparation of this manuscript, the authors used ChatGPT (OpenAI) for limited language editing. The authors take full responsibility for the final content of the manuscript.

\section*{Acknowledgment.} The authors would like to thank Hemanshu Kaul for helpful conversations.  This project was completed by the Coloring Research Group of South Alabama at the University of South Alabama from summer 2025 through fall 2026. The support of the University of South Alabama is gratefully acknowledged.

\bibliographystyle{hplain}
\bibliography{bibliography}{}
\vspace{-5mm}
\end{document}